\documentclass{amsart}
\usepackage{amsfonts,amscd}  

\newtheorem{claim}{}[section]
\newtheorem{theorem}[claim]{Theorem}

\newtheorem{corollary}[claim]{Corollary}

\def\proclaim #1. #2\par{\medbreak
\noindent{\bf#1.\enspace}{\sl#2}\par\medbreak}
\makeatother

\DeclareMathOperator{\ca}{C^*-\text{algebra}}
\DeclareMathOperator{\cas}{C^*-\text{algebras}}

\begin{document}

\title[Complete isometries - an illustration]
{Complete isometries - an illustration of
noncommutative functional analysis}

\date{Aug 23, 2002}   
 
\author{David P. Blecher and Damon M. Hay}
\address{Department of Mathematics, University of Houston, Houston,
TX
77204-3008}
\email[David P. Blecher]{dblecher@math.uh.edu}
\email[Damon Hay]{dhay@math.uh.edu}  
\thanks{This research was supported in part by a grant from
the National Science Foundation}

\begin{abstract} 
This article, addressed to a
general audience of functional analysts, is intended to be an illustration
of a few basic
principles from `noncommutative functional analysis',
more specifically the new field of {\em operator spaces.} 
In our  illustration we show how the classical characterization
of (possibly non-surjective)
isometries between function algebras generalizes to
operator algebras.
We give some variants of
this characterization, and a new proof which has some
advantages.   
\end{abstract}  
 
\maketitle

 
\let\text=\mbox

Intended for Conference Proceedings

\section{Introduction}

The field of operator spaces 
provides a new bridge from 
the world of Banach spaces and function 
spaces, to the world of spaces of operators
on a Hilbert space.   
For researchers in the new field, the philosophical starting point 
is the combination of the following two 
obvious facts.   Firstly, by the Hahn-Banach theorem 
any Banach space $X$ is canonically
linearly isometric to a closed linear subspace 
of $C(K)$, where $K$ is the compact space Ball$(X^*)$.  
Secondly, $C(K)$ is a commutative $C^*$-algebra.
Thus one defines a {\em noncommutative 
Banach space}, or {\em operator space}, to be a 
closed linear subspace $X$ of a possibly noncommutative $C^*$-algebra
$A$.  This simplistic idea becomes much 
more substantive with the addition of some additional 
metric structure.   The point is that if $A$ is any $C^*$-algebra,
then the $*$-algebra $M_n(A)$ of $n \times n$ matrices
with entries in $A$ has a unique norm $\Vert \cdot \Vert_n$
making it a $C^*$-algebra (this follows from the 
well known unicity of $C^*$-norms on a $*$-algebra).  
If $X \subset A$
then $M_n(X)$ inherits this norm $\Vert \cdot \Vert_n$, and 
more precisely we think of an operator space as the
pair $(X,\{ \Vert \cdot \Vert_n \}_n)$.   We usually insist that 
maps between operator spaces are {\em completely bounded},
 where the adjective
`completely' means that we are applying our maps to
matrices too.  Thus if $T : X \rightarrow Y$, then $T$
 is {\em completely contractive} if $T_n$ is 
contractive for all $n \in \mathbb{N}$, where 
$T_n$ is the map $[x_{ij}] \mapsto [T(x_{ij})]$.    
Similarly $T$ is {\em completely isometric}
if $\Vert [T(x_{ij})] \Vert = \Vert [x_{ij}] \Vert$
for all $n \in \mathbb{N}$ and $[x_{ij}] \in M_n(X)$.  
It is an easy exercise (using one of the common expressions
for the 
operator norm of a matrix in $M_n = M_n(\mathbb{C})$)   to 
prove that a linear map $T : X \to Y$ between subspaces of 
$C(K)$ spaces is completely contractive if and only if it is
contractive.  Consequently
such a $T$ is isometric if and only if it is
completely isometric.  

The identification of 
the term `noncommutative Banach space'
with `operator space' may be thought of  
as a relatively recent entry in 
the well known `dictionary' translating terms between
the `commutative' and `noncommutative' worlds.  We spend a paragraph
describing some other entries in this dictionary.
Although these items are for the most part
well known to the point of being tedious, it will be helpful
to collect them here for the dual purpose of 
establishing notation, and for ease of reference later in the paper.
The most well known item 
is of course the fact that the noncommutative
version of a $C(K)$ space is a unital $C^*$-algebra $B$.  
The noncommutative version of a unimodular function in $C(K)$ is
a {\em unitary} $u \in B$ (i.e. $u^* u = u u^* = 1$).
The noncommutative version of a function algebra $A \subset C(K)$ 
`containing constant functions' is a
closed  subalgebra $A$ of a
$C^*$-algebra $B$, with $1_B \in A$.  We call such $A$ a
{\em unital operator algebra}.
For a unital subset ${\mathcal S}$ of a 
$C^*$-algebra $B$, we will take as
a simple noncommutative version of the assertion 
`${\mathcal S} \subset C(K)$ separates points of $K$',
the assertion `the $C^*$-subalgebra of $B$ generated by ${\mathcal S}$ 
(namely, the 
smallest $C^*$-subalgebra of $B$ containing ${\mathcal S}$) equals $B$'.
The analogue of a closed subset $E$ of a compact set $K$
is a quotient $B/I$, where $I$ is a closed two-sided ideal 
in a unital $C^*$-algebra $B$.   
More generally,
unital $*$-homomorphisms $\pi$  between unital  $C^*$-algebras
are the noncommutative version of  continuous functions $\tau$
between compact spaces.  Indeed clearly any such $\tau
: K_1 \to K_2$ gives rise 
to the unital $*$-homomorphism $C(K_2) \to C(K_1)$ 
of `composition with $\tau$', and conversely it is not much harder to 
see that any unital $*$-homomorphism $C(K_2) \to C(K_1)$
comes from a continuous $\tau$ in this way.
Moreover such $\pi$ is 1-1 (resp. onto)
if and only if the corresponding $\tau$ is onto (resp. 1-1).    
Thus the noncommutative version of a
homeomorphism between compact spaces is a
(surjective 1-1) $*$-{\em isomorphism}
between unital $C^*$-algebras.
Coming back to `noncommutative functional analysis',
it is convenient for some purposes
(but admittedly not for others) to view `complete isometries' as the 
noncommutative version of isometries.  It is very important in 
what follows that a 1-1 $*$-homomorphism $\pi :
A \to B$ between $C^*$-algebras,
is by a simple and well known spectral theory argument,
automatically an isometry, and consequently 
(by the same principle applied to $\pi_n$), 
a complete isometry.   Similarly, a $*$-homomorphism 
$\pi : A \to B$ (which is not a priori assumed continuous) is automatically 
completely contractive, and has a closed range which is
a $C^*$-algebra $*$-isomorphic to the   $C^*$-algebra
quotient of $A$ by the obvious two-sided ideal,
namely the kernel of the $*$-homomorphism.
 
The entries we have just described in this `dictionary' 
are all easily justified by well known theorems (for example
Gelfand's characterization of commutative $C^*$-algebras).
That is, if one applies the noncommutative definition
in the commutative world, one recovers exactly the classical object.
Similarly 
one sometimes finds oneself in the very nice
`ideal situation' where one can prove a
theorem or establish a theory in the noncommutative world (i.e.
about operator spaces or operator algebras), which when one 
applies the theorem/theory to objects which
are Banach spaces or function algebras,
one recovers exactly the classical theorem/theory.
An illustration of this point is the Banach-Stone 
theorem.    The following
is a much simpler form of Kadison's characterization 
of isometries between $C^*$-algebras \cite{Kad}:

\begin{theorem} \label{sK}  {\rm (Folklore)} \
 A surjective linear map $T :  A \rightarrow B$ between unital 
$C^*$-algebras is a complete isometry if and only if 
$T = u \pi(\cdot)$, for a unitary $u \in B$ and a 
$*$-isomorphism $\pi : A \to B$.
\end{theorem}  

\begin{proof}  (Sketch.)  
The easy direction  is essentially just 
the fact mentioned earlier that 1-1 $*$-homomorphisms are
completely isometric.  The other direction can be proved  
by first showing (as with 
Kadison's theorem) that $T(1)$ is unitary, so that without loss
of generality $T(1) = 1$.  The well known 
Stinespring theorem has as a simple consequence 
the Kadison-Schwarz inequality $T(a)^* T(a) \leq T(a^* a)$.
Applying this to $T^{-1}$ too yields 
$T(a)^* T(a) = T(a^* a)$, and now the 
result follows immediately from the `polarization 
identity' $a^* b = \frac{1}{4} \sum_{k=0}^3 (a + i^k b)^* (a + i^k b)$.
\end{proof}  
 
Note that if one takes $A = C(K_1)$ and $B = C(K_2)$ 
in Theorem \ref{sK}, and consults the `dictionary' above,
then one recovers exactly the classical 
Banach-Stone theorem.  Indeed as we remarked earlier,
in this case complete isometries are the same thing as
isometries, 
unitaries are unimodular functions, and a
$*$-isomorphism is induced by a homeomorphism between the
underlying compact spaces.   

Indeed consider the following generalization of the 
Banach-Stone theorem: 

\begin{theorem} \label{mat} \cite{Hol,Nov,AF,Math}  Let 
$\Omega$ be compact and Hausdorff, and $A$ a unital function 
algebra.  A linear contraction
 $T :  A \rightarrow C(\Omega)$ is an isometry if and only if 
there exists a closed subset $E$ of $\Omega$,
and two continuous functions $\gamma : E \rightarrow \mathbb{T}$
and $\varphi : E  \rightarrow \partial A$, with $\varphi$
surjective, such that for all $y \in E$
 $$T(f)(y) \; = \; \gamma(y) f(\varphi(y)).$$ \end{theorem}

Here $\partial A$ is the Shilov boundary of $A$ (see 
Section 2).   We have supposed that $T$ maps into a
`selfadjoint function algebra' $C(\Omega)$; however since any
function  algebra is a unital subalgebra of a `selfadjoint' 
one, the theorem also applies to 
isometries between unital function algebras.
If $A$ is 
a $C(K)$ space too,  then $\partial A = K$ and then the 
theorem above is called {\em Holsztynski's theorem}.    We refer the 
reader to \cite{JP} for a survey of such variants on the 
classical Banach-Stone theorem.      

Often the transition from the `classical' to the `noncommutative' involves
the introduction of much more algebra.  Next we appeal to our dictionary 
above to give an equivalent restatement of Theorem 
\ref{mat} in more algebraic 
language. 

\begin{theorem} \label{matr} {\rm (Restatement of Theorem \ref{mat})} \
 Let $A, B$ be unital function algebras, with $B$ selfadjoint.
A linear contraction  $T : A \rightarrow B$ is an isometry if and only if
\begin{itemize}
\item [(A)]  
there exists a closed ideal $I$ of $B$, a unitary $u$ in the quotient
$C^*$-algebra $B/I$, and a unital 1-1 $*$-homomorphism $\pi : A \to B/I$,
such that $q_I(T(a)) = u \pi(a)$ for all $a \in A$. 
\end{itemize}     
Here $q_I$ is the canonical quotient $*$-homomorphism $B \to B/I$.     
\end{theorem}

In light of Theorems  \ref{sK} and \ref{matr}
one would imagine that for any complete isometry 
$T : A \to B$ between unital operator   algebras, 
the condition (A) above 
should hold verbatim.  This would give a pretty noncommutative 
generalization of Theorem \ref{matr}.    Indeed if Ran $T$
is also a unital operator algebra, then this is true
(see eg. B.1 in \cite{BSh}).
However, it is quite easily seen
that such a result cannot hold generally. 
For example, let $M_n = M_n(\mathbb{C})$;
 for any $x \in M_n$ of norm $1$,
the map $\lambda \mapsto \lambda x$ is a complete isometry from
$\mathbb{C}$ into $M_n$.   Now $M_n$ is simple (i.e. has no
nontrivial two-sided ideals), and so if the result above was valid then it
follows immediately that $x = u$.    This is obviously not satisfactory.

To resolve the dilemma presented in the last paragraph, we have offered in
\cite{BlecherHay} several alternatives.    For example, one may replace
the quotient $B/I$ by a quotient of a certain $*$-subalgebra of $B$.
The desired relation $q_I(T(a)) = u \pi(a)$ then requires $u$ to 
be a unitary in a certain $C^*$-triple system (by which
we mean a subspace $X$ of a $C^*$-algebra $A$ with $X X^* X \subset X$).
Or, one 
may replace the quotient $B/I$ by a quotient $B/(J + J^*)$, where
$J$ is a one-sided ideal of $B$.   Such a
quotient is not an algebra, but is an `operator system'
(such spaces have been important in the deep work of Kirchberg
(see \cite{Ki,KW} and references therein).   Alternatively,
one may replace such quotients altogether,
with certain subspaces of the second dual
$B^{**}$ defined in terms of certain orthogonal projections of
`topological significance' (i.e. correspond to 
characteristic functions of closed sets in $K$ if
$B = C(K)$) in the second dual
$B^{**}$ (which is a  von Neumann algebra \cite{Ped}).
The key point of all these arguments, and indeed a key
approach to Banach-Stone theorems for linear maps between
function algebras, $C^*$-algebras or operator algebras,
is the basic theory of $C^*$-triple systems and triple morphisms,
and the basic properties
of the {\em noncommutative Shilov boundary}  or {\em triple envelope}
of an operator space.   These important and beautiful ideas
originate in the work of
Arveson, Choi and Effros, Hamana, Harris, Kadison,
Kirchberg, Paulsen, Ruan, and others.
Indeed our talk at the conference spelled out
these ideas and their connection with 
the Banach-Stone theorem; and the background ideas
are developed at length in a  book the first author
is currently writing with Christian Le Merdy 
\cite{BLM} (although 
we do not characterize non-surjective complete isometries there).  
Moreover, a description of our work from this perspective,
together with many related results, may
be found in \cite{IOBS}.
Thus we will content ourselves here with a survey of some related 
and interesting topics, and with a new and
self-contained proof
of some characterizations of complete isometries between unital 
operator algebras which do not 
appear elsewhere.   This proof has several advantages,
for example the projections arising naturally with this 
approach seem to be more useful for some purposes.
Also it will allow us to avoid any explicit mention of the theory
of triple systems (although this is playing a silent role nonetheless).

We also
show how such noncommutative results are generalizations of the older
characterizations of into isometries between 
function algebras or $C(K)$ spaces.   We thank A. Matheson 
for telling us about these results.  In the final section we  present some
evidence towards the claim that (general) isometries between 
operator algebras are not the correct
noncommutative generalization of isometries between function 
algebras. 

For the reader who wants to learn more 
operator space theory we have listed
some general texts in our bibliography. 

\section{The noncommutative Shilov boundary}

At the present time the appropriate `extreme point' theory is not 
sufficiently developed to be extensively used 
in noncommutative functional analysis.  Although several major
and beautiful pieces are now in place,
 this is perhaps one 
of the most urgent needs in the subject.   However there are  
good substitutes for `extreme point' arguments.    One such is
the {\em noncommutative Shilov boundary} of an operator space.    
Recall that if $X$ is a closed subspace of $C(K)$ containing 
the identity function $1_K$ on $K$ and separating points of $K$, 
then the classical Shilov boundary may be defined to be 
the smallest closed subset $E$ of $K$ such that all functions
$f \in X$ attain their norm, or equivalently such that the restriction
map $f \mapsto f_{|_E}$ on $X$ is an isometry.   This boundary is often defined 
independently of $K$, for example if $A$ is a unital function 
algebra then we may define the Shilov boundary as we just did,
but with $K$ replaced by the maximal ideal space of $A$.   In fact we 
prefer to think of the classical Shilov boundary of $X$ 
as a pair $(\partial X,i)$ consisting of an abstract compact 
Hausdorff space $\partial X$, together with an isometry 
$j : X \rightarrow C(\partial X)$ such that $j(1_K) = 1_{\partial X}$
and such that $j(X)$ separates points of $\partial X$,  with the 
following universal property:  
For any other pair $(\Omega,i)$ consisting of a compact
Hausdorff space $\Omega$  
and a complete isometry $i : X \rightarrow C(\Omega)$ which is
unital (i.e. $i(1_K) = 1_A$),
and  such that $i(X)$ separates points of $\Omega$,
there exists a (necessarily
unique) continuous injection $\tau : \partial X \to \Omega$ 
such that $i(x)(\tau(w)) = j(x)(w)$ for all 
$x \in X, w \in \partial X$.
Such a pair $(\partial X,i)$ is easily seen to 
be unique up to an appropriate 
homeomorphism.   The fact that such $\partial X$ exists is the 
difficult part, and proofs may be found in books
on function algebras (using extreme point arguments).
  
Consulting our `noncommutative dictionary' in Section 1, and thinking 
a little about the various correspondences there, it will 
be seen that the noncommutative version of this
 universal property above should read as follows.  Or at any rate,
the following noncommutative statements, when applied to 
a unital subspace $X \subset C(K)$, will imply the 
universal property of the classical Shilov boundary discussed 
above.    Firstly,
a {\em unital operator space} is a pair $(X,e)$ consisting
of an operator space $X$ with fixed element $e \in X$,
such that there exists a linear complete isometry 
$\kappa$ from $X$ into a unital $\ca$ $C$ with $\kappa(e) = 1_C$.
A `noncommutative Shilov boundary' would correspond to
a pair $(B,j)$ consisting of a unital $\ca$ $B$
and a complete isometry $j : X \rightarrow
B$ with $j(e) = 1_B$, and 
whose range generates $B$ as a $\ca$, with the following
universal property:  For any other pair $(A,i)$ consisting of a unital $\ca$
and a complete isometry $i : X \rightarrow A$ which is 
unital (i.e. $i(e) = 1_A$),
and whose range generates $A$ as a $\ca$, there exists a (necessarily
unique, unital, and surjective) *-homomorphism $\pi :
A \rightarrow B$ such that $\pi \circ i = j$.
Happily, this turns out to be true.  The 
existence for any unital operator space $(X,e)$
of a pair $(B,j)$ with the universal property above is 
of course a theorem, which  
 we call the Arveson-Hamana theorem \cite{Arv,Ham1} (see
\cite{BSh} for complete details).
As is customary we write $C^*_e(X)$ for $B$ or $(B,j)$, this is
the `$C^*$-envelope of $X$'.   It is essentially
unique, by the universal property.    
If $X = A$ is a unital operator algebra (see Section 1
for the definition of this),
then $j$ above
is forced to be a
homomorphism (to see this, choose an $i$ which is a 
homomorphism, and use the universal 
property). Thus  $A$ may be considered a unital subalgebra of
 $C^*_e(A)$.  If $A$ is already a unital $\ca$, then of course we 
can take $C^*_e(A) = A$.

To help the reader get a little more 
comfortable with these concepts,
 we compute the `noncommutative Shilov boundary' 
in a few simple examples.

\medskip

{\bf Example 1.}  Let $T_n$ be the upper triangular
$n \times n$ matrices.  This is a unital subspace of $M_n$,
and no proper $*$-subalgebra of $M_n$ contains $T_n$.
Let $(B,j)$ be the $C^*$-envelope of $T_n$.  By the 
universal property of the $C^*$-envelope, there is a
surjective $*$-homomorphism $\pi : M_n \to B$ such that 
$\pi(a) = j(a)$ for $a \in T_n$.  The kernel of 
$\pi$ is a two-sided ideal of $M_n$.   However 
$M_n$ has no nontrivial two-sided ideals.  Hence  
$\pi$ is 1-1, and is consequently a $*$-isomorphism,
and we can thus identify $M_n$ with $B$.
Thus $M_n$ is a $C^*$-envelope of $T_n$.  

\medskip

{\bf Example 2.}   Consider the linear
subspace $X$ of $M_3$ with zeroes in the
1-3, 2-3, 2-1, 3-1 and 3-2 entries, and with arbitrary entries 
elsewhere except for the 3-3 entry, which is the average 
of the 1-1 and 2-2 entries.   It is easy to see that
the $C^*$-algebra generated by $X$ inside $M_3$
is $M_2 \oplus {\mathbb C}$.  However this is
not the $C^*$-envelope.  Indeed it is easy to see that 
the 3-3 entry here is redundant, since the norm of 
$x \in X$ is the norm of the upper left $2 \times 2$ block 
of $x$.   This observation can be expanded to show that 
the canonical projection map $M_2 \oplus {\mathbb C} \to M_2$
when restricted to $X$ is a unital complete isometry
from $X$ onto $T_2$ (see Example 1). 
This is the same as saying that if one takes the quotient of 
$M_2 \oplus {\mathbb C}$ by its ideal $0_2 \oplus
{\mathbb C}$, then one obtains $M_2$, which by Example 1
is the $C^*$-envelope.  

Indeed this is typical when calculating the $C^*$-envelope of a
 unital subspace $X$ of $M_n$.  The  $C^*$-algebra generated by $X$
is a finite dimensional unital $C^*$-algebra.  However such
$C^*$-algebras are all $*$-isomorphic to a finite 
direct sum $B$ of full `matrix blocks' $M_{n_k}$.  Some of these
blocks are redundant.  That is, if $p$ is the central projection 
in $B$ corresponding to the identity matrix of this block, 
then $x \mapsto x(1_B - p)$ is completely isometric.    
If one eliminates such blocks then the 
remaining direct sum of blocks is the $C^*$-envelope.
 
\medskip
 
{\bf Example 3.}  Let $B$ be
a unital $C^*$-algebra.
Consider the unital subspace ${\mathcal S}(B)$ 
of the $C^*$-algebra $M_2(B)$ consisting of matrices
$$\left[ \begin{array}{ccl} \lambda 1 & x \\ y^* & \mu 1 \end{array}
\right]$$
for all $x, y \in B$ and $\lambda, \mu$ complex scalars.
We claim that $M_2(B)$ is the $C^*$-envelope $C$
of ${\mathcal S}(B)$,
and we will prove this using a similar idea to 
Example 1 above.  Namely, first note that 
$M_2(B)$ has no proper $C^*$-subalgebra containing 
${\mathcal S}(B)$,  Thus by the Arveson-Hamana theorem there 
exists a $*$-homomorphism $\pi : M_2(B) \to 
C$ which possesses a property which we will 
not repeat, except to say that it certainly ensures that 
$\pi$ applied to a matrix with zero entries except for 
a nonzero entry in the 1-2 position, is nonzero.
It suffices as in Example 1 to show that
Ker $\pi = \{ 0 \}$.   Suppose that $\pi(x) = 0$ for 
a $2 \times 2$ matrix $x \in M_2(B)$.   Let $E_{ij}$ be the four
canonical basis matrices for $M_2$, thought 
of as inside $M_2(B)$.  
Then $\pi(E_{1i} x E_{j2}) = 
\pi(E_{1i}) \pi(x) \pi(E_{j2}) = 0$ for $i,j = 1,2$.
Thus by the fact mentioned 
above about the 1-2 position, we must have $E_{1i} x E_{j2} = 0$.  
Thus $x = 0$.

\medskip
 
In fact a variant of the $C^*$-envelope or `noncommutative
Shilov boundary' can be defined for any operator space $X$.
This is the {\em triple envelope} of Hamana (see \cite{Ham2}).
This is explained in much greater detail in \cite{BSh}, together 
with many applications.   For example it is intimately connected 
to the `noncommutative $M$-ideals' recently introduced in \cite{BEZ}.  
This `noncommutative Shilov boundary' is, as we mentioned in
Section 1, a key tool for proving various Banach-Stone type theorems. 
However in the present article we shall only need 
the variant described earlier in this section.
 
\section{Complete isometries between operator algebras}
  
We begin this section with  
a collection of very well known and simple
facts about closed two-sided ideals
$I$ in a $C^*$-algebra $A$, and about the 
quotient $C^*$-algebra $A/I$.   We have 
that  $I^{\perp\perp}$ is a weak* closed 
two-sided ideal in the von Neumann algebra $A^{**}$, and 
there exists a unique orthogonal  projection $e$ in the 
center of $A^{**}$  with $I^{\perp\perp} = A^{**} (1-e)$.  
 The projection $1-e$ is called the {\em support projection } for
$I$, and $1-e$  may be taken to be
the weak* limit in $A^{**}$ of any contractive
approximate identity for $I$.  Thus it follows that
$A^{**}/I^{\perp\perp} \cong A^{**}e$ as $C^*$-algebras.  Therefore
also  $$A/I \subset (A/I)^{**} \cong A^{**}/I^{\perp\perp} \cong A^{**} e$$
as $C^*$-algebras.  Explicitly, the composition of all these 
identifications is a 1-1 $*$-homomorphism taking an $a + I$ in $A/I$, to
$\hat{a} e = e \hat{a} e$ in $A^{**}$.   Here $\; \hat{} \; \; $
is the canonical embedding $A \to A^{**}$ (which we 
will sometimes suppress mention of).  
Thus $A/I$ may be regarded as
a $C^*$-subalgebra of $A^{**}$, or of the $C^*$-algebra $eA^{**}e$.

We next illustrate the main idea of our theorem
with a simple special case.
(The following appeared as part of
Corollary 3.2 in the 
original version of \cite{BlecherHay}, with the 
proof left as an exercise).
Suppose that $T : A \to B$ is a 
complete isometry between unital $C^*$algebras,
and suppose that $T$ is unital too, that is $T(1) = 1$.
Let $C$ be the $C^*$-subalgebra of $B$ generated by $T(A)$.
Applying the Arveson-Hamana theorem\footnote{We remark in
passing that one does not need the full strength of the Arveson-Hamana theorem
here, one may use the much simpler 
\cite{CE} Theorem 4.1.}
 we obtain a surjective $*$-homomorphism
$\theta : C \to A$ such that $\theta(T(a)) = a$ for all $a \in A$.  
If $I$ is the kernel of 
the mapping $\theta$, then $C/I$ is a unital $C^*$-algebra $*$-isomorphic
to $A$.  Indeed there is the canonical
$*$-isomorphism $\gamma : 
A \to C/I$ induced by $\theta$, taking $a$ to $T(a) + I$.
The next point is that $C/I$ may be viewed as we mentioned 
a few paragraphs back, as a $C^*$-subalgebra of $C^{**}$,
and therefore also of $B^{**}$.  
Indeed if $e$ is the central projection in $C^{**}$ mentioned there,
then $C/I$  may be viewed as a $C^*$-subalgebra of
$e C^{**} e \subset e B^{**} e \subset B^{**}$.
In view of the last fact,
the map $\gamma$ induces an 1-1 $*$-homomorphism
$\pi : A \to B^{**}$ taking an element $a \in A$ to the element of $B^{**}$
which equals  
\begin{equation} \label{eq1}
\widehat{T(a)} e \; = \; e \widehat{T(a)} 
 \; = \; e \widehat{T(a)} e 
\end{equation} 
(these are equal because $e$ is central in $C^{**}$).  
Conversely,
if $T : A \to B$ is a complete contraction for which there
exists a projection $e \in B^{**}$ such that 
$e \widehat{T(a)} e$ is
a  1-1 $*$-homomorphism
$\pi$, 
then for all $a \in A$,
$$\Vert T(a) \Vert \geq \Vert e \widehat{T(a)} e \Vert
= \Vert \pi(a) \Vert = \Vert a \Vert$$
using the fact mentioned earlier that 1-1 $*$-homomorphisms
are necessarily isometric.   
Thus $T$ is an isometry, and a similar argument shows that
it is a complete isometry.     Thus we have characterized
unital complete isometries $T : A \to B$.  

If $H$ is a Hilbert space on which we have represented the
von Neumann algebra $B^{**}$ as a weak* closed unital
$*$-subalgebra, then $B$ may be viewed also as a
unital $C^*$-subalgebra of $B(H)$, whose weak* closure in
$B(H)$ is (the copy of) $B^{**}$.   In this case we shall
say that $B$ is represented on $H$ {\em universally}.  
(The explanation for this term is that the well-known 
`universal representation' $\pi_u$ of a $C^*$-algebra is
`universal' in our sense, and conversely 
if  $\pi$ is a
representation which is `universal' in our sense
then $\pi(B)''$ is isomorphic to $\pi_u(B)'' \cong B^{**}$.
See \cite{Rief} Section 1.)  
If, further, $e \in B^{**}$  is a projection for which 
(\ref{eq1}) holds, then with respect to the 
splitting $H = eH \oplus (1-e)H$  we may write
$$T(a) = \left[ \begin{array}{ccl}
\pi(\cdot) & 0 \\ 0 & S(\cdot) \end{array} \right],$$
for all $a \in A$.     We will see that this is essentially
true even if $T(1_A) \neq 1_B$:  
   
\begin{theorem} \label{BH}  
Let $T : A \rightarrow B$ be a completely contractive linear map
from a unital operator algebra into
a unital $\ca$.  Then the following are equivalent:
\begin{itemize}
\item [(i)]   $T$ is a complete isometry,
\item [(ii)] There is a
partial isometry $u \in B^{**}$ with initial projection $e \in B^{**}$, 
and a (completely isometric) 1-1 
$*$-homomorphism $\pi : C^*_e(A) \to e B^{**} e$ with 
$\pi(1) = e$, such that for all $a \in A$ 
$$\widehat{T(a)} e = u \pi(a) \; \; \;
\text{and}  \; \; \;  \pi(a) = u^* \widehat{T(a)}.$$
Moreover $e$ may be taken to be 
a `closed projection' (see
\cite{Ped} 3.11, and the discussion towards the
end of our proof).
\item [(iii)]  If $H$ is a Hilbert space 
on which $B$ is represented universally,
then there exist two closed subspaces $E, F$ of the Hilbert
space $H$, a  1-1     $*$-homomorphism $\pi :  C^*_e(A) \to B(E)
$ with $\pi(1) = I_E$, and a unitary $u : E \to F$, such that
$$T(a)_{|_E} = u \pi(a),$$
and $T(a)_{|_{E^\perp}} \subset F^\perp$, for all $a \in A$.
Here $E^\perp$ for example is the orthocomplement of
$E$ in $H$.
  \item [(iv)]  If $H$ is as in (iii),
then there exists two closed subspaces $E, F$ of $H$,  
a unital 1-1 $*$-homomorphism $\pi : C^*_e(A) \to B(E)$,
a complete contraction $S :  C^*_e(A) \to B(E^\perp,F^\perp)$, and 
unitary operators $U : E \oplus F^\perp \to H$ and $V : H 
\to E \oplus E^\perp$,   
such that 
$$
T(a)
\; = \; U \left[ \begin{array}{ccl} \pi(a) & 0 \\
0 & S(a) \end{array} \right] V$$
for all $a \in A$.
\item [(v)]   There is a left ideal $J$ of $B$, a
1-1 *-homomorphism $\pi$ from $C^*_e(A)$ into a
unital subspace of $B/(J + J^*)$ which is a $C^*$-algebra,
and a `partial isometry' $u$ in $B/J$
such that $$q_J(T(a)) = u \pi(a)
\; \; \; \; \; \; \; \; \& \; \; \;
 \; \; \; \; \; \pi(a)  = u^* q_J(T(a))$$
for all $a \in A$, where $q_J$ is the canonical quotient map
$B \rightarrow B/J$.
\end{itemize} 
\end{theorem} 

Before we prove the theorem,
we make several remarks.  First, we have taken $B$ to 
be a $C^*$-algebra; however since any 
unital operator algebra is a
unital subalgebra of a unital $C^*$-algebra this is
not a severe restriction.
We also remark that there are several 
other items that one might add to such a list of 
equivalent conditions.   See \cite{BlecherHay,BL}.
Items (ii)-(iv),
and the proof given below of their equivalence
with (i), are new.  
We acknowledge that
we have benefitted from a suggestion that we use the 
Paulsen system to prove the result.   
This approach is an obvious one to those working in this
area (Ruan and Hamana used a variant of it in their
work in the '80's on complete isometries and triple morphisms
\cite{Ru,Ham2}).  However  we had not pushed through this approach 
in the original 
version of \cite{BlecherHay} because this method does 
not give several of the results there as immediately.
Statement (v) above has been simply copied from 
\cite{BlecherHay,BL} without proof or explanation.
We have listed it here simply because Theorem \ref{matr}
may be particularly easily derived from it as the special case when
$A$ and $B$ are commutative (see comments below).
Note that (iii) above resembles  Theorem \ref{mat}  superficially.

\begin{proof}   The fact that the other conditions
all imply (i) is easy, following the idea in the paragraph 
above the theorem, namely by using the fact that a 
1-1 $*$-homomorphism is completely isometric.     

In the remainder of the proof we
suppose that $T$ is a complete isometry.
We view $A$ as a unital subalgebra of $C^*_e(A)$ as
outlined in Section 3.  We 
define a subset ${\mathcal S}(B)$ of $M_2(B)$ as in 
Example 3 in Section 2.  Similarly define a subset
${\mathcal S}(T(A))$ of ${\mathcal S}(B)$ 
using a similar formula (note that 
${\mathcal S}(T(A))$ has 1-2 entries taken from
$T(A)$ and 2-1 entries taken from $T(A)^*$).
Similarly we define 
the subset ${\mathcal S}(A)$ of the 
$C^*$-algebra $M_2(C^*_e(A))$
(i.e. ${\mathcal S}(A)$ has scalar diagonal 
entries and off diagonal entries from $A$ and $A^*$).
We write $1 \oplus 0$ for the matrix in ${\mathcal S}(A)$
with $1$ as the 1-2 entry and zeroes elsewhere.
Similarly for  $0 \oplus 1$.
We also use these expressions for the 
analogous matrices in ${\mathcal S}(B)$.
The map $\Phi : {\mathcal S}(A) \to 
{\mathcal S}(T(A)) \subset M_2(B)$
taking 
$$\left[ \begin{array}{ccl} \lambda 1 & x \\ y^* & \mu 1 \end{array}
\right] \; \mapsto \; \left[ \begin{array}{ccl} \lambda 1 & T(x)
 \\ T(y)^* & \mu 1 \end{array}
\right]$$
is well known to be a unital complete isometry (this is the well known
Paulsen lemma, see the proof of
7.1 in \cite{P}).   Let $C$ be the 
$C^*$-subalgebra of $M_2(B)$ generated by ${\mathcal S}(T(A))$.
The $C^*$-envelope of ${\mathcal S}(A)$ is well known 
to be $M_2(C^*_e(A))$ (see Example 3 in Section 2 where we 
proved this in the case that $A$ is already a 
$C^*$-algebra, or for example \cite{BSh} Proposition
4.3 or \cite{Zh}).  Thus by the 
Arveson-Hamana theorem we obtain a surjective $*$-homomorphism
$\theta : C \to M_2(C^*_e(A))$ such that $\theta \circ \Phi$ is simply the 
canonical embedding of ${\mathcal S}(A)$ into $M_2(C^*_e(A))$.
As in the special case
considered above the theorem, we let $I_0$ be the kernel of
the mapping $\theta$, then $C/I_0$ is a unital $C^*$-algebra $*$-isomorphic
to $M_2(C^*_e(A))$.  Indeed there is the canonical
$*$-isomorphism $\gamma : M_2(C^*_e(A)) \to C/I_0$ induced by $\theta$, taking 
$$\left[ \begin{array}{ccl} \lambda 1 & x \\ y^* & \mu 1 \end{array}
\right] \; \mapsto \; \left[ \begin{array}{ccl} \lambda 1 & T(x)
 \\ T(y)^* & \mu 1 \end{array}
\right] \; + \; I_0.$$              
As in the simple case above the 
theorem,  $C/I_0$ may be viewed as a $C^*$-subalgebra of 
$p_0 C^{**} p_0$, for a central projection $p_0 \in C^{**}$
(namely, the complementary projection to
the support projection of $I_0$).  
Now $p_0 C^{**} p_0 \subset C^{**} \subset M_2(B)^{**}$,
and it is well known that
$M_2(B)^{**} \cong M_2(B^{**})$ as $C^*$-algebras.
Thus we may think of $C^{**}$ as a $C^*$-subalgebra of
$M_2(B^{**})$.  Also,
$C^{**}$ contains $C$ as a $C^*$-subalgebra, and 
the projections $1 \oplus 0$ and $0 \oplus 1$ in $C$ 
correspond to the matching diagonal projections 
$1 \oplus 0$ and $0 \oplus 1$ in $M_2(B^{**})$.
These last projections therefore commute with
$p_0$, since $p_0$ is central in $C^{**}$, which  immediately
implies that 
$p_0$ is a diagonal sum $f \oplus e$ of two  
orthogonal projections $e,  f \in B^{**}$.   
Thus we may write 
the $C^*$-algebra $p_0 M_2(B^{**}) p_0$
as the $C^*$-subalgebra
$$\left[ \begin{array}{ccl} f B^{**} f & f B^{**} e \\
e B^{**} f & e B^{**} e \end{array}
\right]$$
of $M_2(B^{**})$.    
We said above that $C/I_0$ may be regarded as a 
$C^*$-subalgebra of the subalgebra
$p_0 M_2(B^{**}) p_0$ of $M_2(B^{**})$.
Thus the map $\gamma$ induces a  1-1    $*$-homomorphism
$\Psi : M_2(C^*_e(A)) \to M_2(B^{**})$.   
It is easy to check that 
$\Psi(1 \oplus 0) = f \oplus 0$ and $\Psi(0 \oplus 1) =
0 \oplus e$.  
Since $\Psi$ is a $*$-homomorphism it follows that 
$\Psi$ maps each of the four corners of $M_2(C^*_e(A))$ to the 
corresponding corner of $p_0 M_2(B^{**}) p_0
\subset M_2(B^{**})$.  We let $R : C^*_e(A) \to f B^{**} e$ be the
restriction of $\Psi$ to the `1-2-corner'.  Since $\Psi$ is
 1-1, it follows that $R$ is  1-1.
If $\pi$ is the restriction of
$\Psi$ to the `2-2-corner', 
then $\pi$ is a *-homomorphism $C^*_e(A) \to e B^{**} e$ taking
$1_A$ to $e$.  Applying the $*$-homomorphism $\Psi$ 
to  the identity
$$\left[ \begin{array}{ccl}  0 & 0 \\ 1 & 0 \end{array}
\right] \left[ \begin{array}{ccl} 0 & 1 \\ 0 & 0 \end{array}
\right] =
\left[ \begin{array}{ccl} 0 & 0 \\ 0 & 1 \end{array}
\right]
$$
we obtain that $u = R(1)$ is a
partial isometry, with $u^* u = \pi(1) = e$.
Similarly $u u^* = f$.   A similar argument shows that 
$R(a) = R(1) \pi(a)$ for all $a \in C^*_e(A)$.
Thus $u^* R(a) = u^* u \pi(a) = \pi(a)$ for all $a \in 
C^*_e(A)$.  
  
Next, we observe that $\Psi$ takes the matrix $z$ which is zero 
except for an $a$ from $A$ in the 1-2-corner, to the matrix 
$w = p_0 \widehat{\Phi(z)} p_0$.  Since 
$\widehat{\Phi(z)} \in C^{**}$ and $p_0$ is in the center of that 
algebra, we also have $w =  \widehat{\Phi(z)} p_0 =
p_0 \widehat{\Phi(z)}$.  Also $w$ viewed as
a matrix in $M_2(B^{**})$ has zero entries except in 
the 1-2-corner, which (by the last sentence) equals  
$$f \widehat{T(a)} e = \widehat{T(a)} e = f \widehat{T(a)}.$$
Also using these facts and a fact from the end of the last 
paragraph we have
$$u^* \widehat{T(\cdot)} = R(1)^* \widehat{T(\cdot)} =
(f \widehat{T(1)} e)^* \widehat{T(\cdot)}
= e \widehat{T(1)}^* f \widehat{T(\cdot)}
=  e \widehat{T(1)}^* \widehat{T(\cdot)} e 
= u^* R(\cdot) = \pi.$$
Thus 
$$\widehat{T(\cdot)} e = 
f \widehat{T(\cdot)} = u u^* \widehat{T(\cdot)} =
u \pi(\cdot).$$
We have now also established most of (ii).
One may deduce (iii) from (ii) by viewing $B \subset B^{**} \subset B(H)$,
and setting $E = eH$, and $F = (uu^*) H$.
We also need to use facts from the
proof above such as $u^* u = e$.
Clearly (iv) follows from (iii).   As we said above, we will 
not prove (v) here.

Claim: if $e$ is the projection
in (ii) above, then $1-e$
is the support projection for a closed ideal $I$
of a unital $*$-subalgebra $D$ of $B$.   Equivalently
(as stated
at the start of this section), there is a
(positive increasing) contractive approximate identity 
$(b_t)$ for $I$, with $b_t \to 1-e$ in the weak* topology.
This claim shows that $1-e$ is
an `open projection' in $B^{**}$, so that 
$e$ is a closed projection, as will be 
obvious to operator algebraists from \cite{Ped} section 3.11
say.  For our other readers we note that for what comes later 
in our paper, one can replace the assertion
about closed projections in the statement of 
Theorem  \ref{BH} (ii) with the statement in the Claim above.

To prove the Claim, recall from
our proof that $p_0 = f \oplus e = 1_C - p_1$, where
$p_1$ is
the support projection for a closed ideal $I_0$
of $C$.   Thus $p_1 = (1-f) \oplus (1-e)$.
As stated at the start of Section 3,
$p_1$ is the weak* limit in $C^{**}$, and hence also
in $M_2(B^{**})$, of a
contractive approximate identity $(e_t)$ of $I_0$.  By
the separate weak* continuity of the product
in a von Neumann algebra, it follows that
the net
$b_t = (0 \oplus 1) e_t (0 \oplus 1)$
has weak* limit $(0 \oplus 1) p_1 (0 \oplus 1) = 0 \oplus (1-e)$.
Viewing these as expressions in $B$, the above says that
$b_t \rightarrow 1-e$ weak* in $B^{**}$.
View $(0 \oplus 1) C (0 \oplus 1)$ as a $*$-subalgebra $D$ of $B$,
and view $(0 \oplus 1) I_0 (0 \oplus 1)$ as a two sided ideal $I$ in $D$.
It is easy to see that  $(b_t)$ is
a contractive approximate identity of $I$.
Thus it follows that $1-e$ is the support
projection of the ideal $I$.
\end{proof}

Some applications of results such as Theorem \ref{BH} may be found in
\cite{BL}.

Next we discuss briefly
the relation between our {\em noncommutative}
characterization of complete isometries (for example Theorem 
\ref{BH} above), and Theorem \ref{matr}.   Our point is not to 
provide another proof for Theorem \ref{matr} - the best existing proof
is certainly short and elegant.  Rather we simply wish to 
show that the noncommutative result contains \ref{matr}.
Indeed Theorem \ref{matr} quite easily follows
from Theorem \ref{BH} (v).
Since however we did not prove
Theorem \ref{BH} (v), we give an alternate proof.

\begin{corollary}
\label{impc}   Let $A, B$ be a unital function 
algebras, with $B$ selfadjoint.  Then condition (ii) in
Theorem \ref{BH} implies condition (A) in 
Theorem \ref{matr}.
\end{corollary} 

\begin{proof}  
By hypothesis,
$T(\cdot) e = u \pi(\cdot)$, and $u^* u = e = \pi(1)$ so that 
$u = u \pi(1) = T(1) e$.  Thus $e T(1)^* T(\cdot) e
= u^* u \pi(\cdot) = \pi(1) \pi(\cdot) = \pi$, so that
Ran $\; \pi \subset e \hat{B} e = \hat{B} e$ (note
$B^{**}$ is commutative in this case).  
From  \cite{Ped} 3.11.10 for example, the `closed projection' 
$e$ in
$B^{**}$ corresponds to a  closed ideal $J$ in $B$
whose support projection is $1-e$.
Alternatively, to avoid quoting facts from \cite{Ped}, we will
also deduce 
this from the `Claim' towards the end of the proof of 
Theorem \ref{BH}.
If $I$ is the ideal in that Claim,
let $J$ be the closed ideal in $B$ generated by $I$. 
Since $J = B I$, the contractive approximate identity 
of $I$ is a right contractive approximate identity of $J$.
Thus $J$ has support projection $1-e$ too, by the first paragraph of 
Section 3 above.   

By facts in the just quoted paragraph, we have a canonical unital 1-1 map
$\eta : B/J \to B^{**}$
taking the equivalence class $b + J$ of $b \in B$ to
$e b e$. 
Indeed in this commutative case 
we see by inspection that $\eta$ is a $*$-homomorphism
from the $C^*$-algebra $B/J$ onto
the $C^*$-subalgebra $M = eBe$ of $B^{**}$.   Define 
$\theta(a) = \eta^{-1}(\pi(a))$, this 
is a 1-1 *-homomorphism $A \to B/J$.   Since
$\pi(1) = e$,  $\theta$ is a unital map too.       
Since $u u^* = u^* u = 
e,$ 
$u$ is unitary in $M$, and so $\gamma = \eta^{-1}(u)$
is unitary in $B/J$.   
Note also that $T(a) e = \eta(T(a) + J)$.
Applying $\eta^{-1}$ to the equation
$T(\cdot) e = u \pi(\cdot)$, we obtain
$q_J(T(a)) = \gamma \; \theta(a)$, that is, 
 condition (A) in Theorem \ref{matr}.
\end{proof}

If one attempts to use the ideas above
to find a characterization analogous to condition (A) from 
Theorem \ref{matr} but in the noncommutative case, it seems to us that one 
is inevitably led to a condition such as (v) in Theorem \ref{BH}.   

We address a paragraph
to experts, on generalizations of the
 proof of Theorem \ref{BH}.  Consider
a complete isometry between possibly non-unital
$C^*$-algebras.  Or much more generally,
suppose that $T$
is a complete isometry from an operator space $X$ into
a $C^*$-triple system $W$.  One may form the
so called `linking $C^*$-algebra' of $W$, with the
identities of the `left and right algebras of $W$'
adjoined.  Call this ${\mathcal L}'(W)$.   As in the
proof of Theorem \ref{BH} we think of ${\mathcal S}(W)
\subset {\mathcal L}'(W)$.   Similarly, if $Z$ is the
`triple envelope' of $X$ (or if $X = Z$ is already
a $C^*$-algebra or $C^*$-triple system), then we may consider
${\mathcal S}(X) \subset {\mathcal S}(Z) \subset {\mathcal L}'(Z)$.
As in the proof of Theorem
 \ref{BH} we obtain firstly a unital complete isometry
$\Phi : {\mathcal S}(X) \to {\mathcal S}(T(X))
\subset {\mathcal L}'(Z)$, and then
a unital 1-1 $*$-homomorphism  $\pi : {\mathcal L}'(Z)
\to {\mathcal L}'(W)^{**}$.  By looking at the
`corners' of $\pi$ we obtain projections
$e, f$ in certain second dual von Neumann algebras,
so that $f T(\cdot) e$ is (the
restriction to $X$ of a completely isometric) a 1-1
triple morphism into $W^{**}$.  In fact we have precisely
such a result in \cite{BlecherHay} (see Section 2 there),
but the {\em key point} is
that the new proof gives different projections $e, f$,
which are more useful for some purposes.

\section{Complete isometries versus isometries}
 
Finally, as promised we discuss why we believe that in this
setting of nonsurjective maps  between  $\cas$ say,
general isometries are not the `noncommutative analogue'
of isometries between function algebras.
The point is simply this.  In the function algebra
case we can say thanks to Holsztynski's theorem that 
the isometries are essentially the maps composed of two 
disjoint pieces $R$ and $S$,
where $R$ is isometric and `nice', and $S$ is
contractive and irrelevant.  However 
at the present time it looks to us unlikely
that there ever will be such a result valid for 
general nonsurjective isometries between general $\cas$.
The chief evidence we present for this
assertion is the 
very nice complementary work of Chu and Wong
\cite{CW} on isometries (as opposed to complete
isometries) $T : A \to B$ between $C^*$-algebras.
They show that for such $T$ there is a largest projection $p \in B^{**}$
such that $T(\cdot)p$ is some kind of Jordan triple morphism.
This appears to be the correct `structure theorem',
or version of Kadison's theorem \cite{Kad}, for 
nonsurjective  isometries.   However as they show, the `nice piece' 
$R = T(\cdot)p$ is very often trivial (i.e. zero), and is 
thus certainly not isometric.   Thus this approach
is unlikely to ever
yield a {\em characterization} of isometries.    
A good example is $A = M_2$, the smallest
noncommutative $C^*$-algebra.  Simply because
$A$ is a Banach space there exists, as in the discussion in
the first paragraph of our paper, a linear isometry of
$A$ into a $C(K)$ space.   However it is easy to see
that there is no nontrivial
$*$-homomorphism or Jordan homomorphism from $A$ into a
commutative $C^*$-algebra.  Such an isometry is uninteresting, 
because the interesting
`nice part' is zero.  Thus we imagine that 
the `good noncommutative notions of isometry' are either
complete isometries or the closely related class of maps
for which the piece $T(\cdot)p$ from \cite{CW}
is an isometry.   
 
This leads to two questions.  Firstly, can one independently
characterize the last mentioned class?  Secondly,
if $T$ is a complete isometry, then is the projection $p$
in the last paragraph equal
(or closely related) to our projection $e$ above?

\end{document}